\documentclass[12pt]{paper}

\author{Benjamin J. Wilson}
\title{Imaginary Highest-Weight Representation Theory and Symmetric Functions}

\theoremstyle{definition}

\newtheorem*{theoremnonum}{Theorem}
\newtheorem*{conjecturenonum}{Conjecture}
\newtheorem{theorem}[equation]{Theorem}
\newtheorem{propn}[equation]{Proposition}
\newtheorem{lemma}[equation]{Lemma}
\newtheorem{corollary}[equation]{Corollary}
\newtheorem{remark}[equation]{Remark}

\newtheorem{conjecture}[equation]{Conjecture}

\newcommand{\K}{\Bbbk}
\newcommand{\FieldC}{\mathbb{C}}
\newcommand{\FieldF}{\mathbb{F}}
\newcommand{\Z}{\mathbb{Z}}

\newcommand{\End}[1]{\text{End}\hspace{0.1em}#1}
\newcommand{\g}{\mathfrak{g}}

\newcommand{\NewTerm}[1]{\textit{#1}}
\newcommand{\Note}[1]{\text{\small{\ (#1)}}}
\newcommand{\ENote}[1]{\text{\small{\quad (#1)}}}

\newcommand{\Inverse}[1]{#1^{-1}}
\newcommand{\Span}[1]{\text{span} \set{#1}}
\newcommand{\SpanField}[2]{\text{span}_{#1} \set{#2}}

\newcommand{\Image}{{\text{im\hspace{0.15em}}}}

\newcommand{\UEA}[1]{\text{U}(#1)}

\newcommand{\Dual}[1]{#1^{\ast}}

\newcommand{\Brac}[2]{\lbrack \hspace{0.1em} #1 , #2 \hspace{0.1em} \rbrack} 
\newcommand{\Adjoint}[1]{\text{ad\hspace{0.1em}}#1}
\newcommand{\CSA}{\mathfrak{h}}
\newcommand{\RootSystem}{\Delta} 
\newcommand{\SLTwo}{\text{sl}(2)}
\newcommand{\Verma}[1]{\mathbf{V}(#1)}
\newcommand{\NPlus}{\mathfrak{N}_{+}} 
\newcommand{\NMinus}{\mathfrak{N}_{-}} 
\newcommand{\VermaQuotient}[1]{\mathbf{M}(#1)} 

\newcommand{\Gen}{\mathrm{u_0}} 
\newcommand{\FProd}[2]{\prod_{1 \leqslant i \leqslant #1}{\C{\slf}{#2_i}}} 
\newcommand{\Heis}{\mathscr{H}} 
\newcommand{\Level}[1]{\mathbf{M}_{#1}} 
\newcommand{\F}{\mathscr{F}} 
\newcommand{\E}{\mathscr{E}} 
\newcommand{\SymGroup}[1]{{\text{Sym}(#1)}} 
\newcommand{\PowerSum}[1]{\mathbf{p}({#1})} 
\newcommand{\ElemSymFunc}[1]{\upvarepsilon_{#1}} 
\newcommand{\Sign}[1]{\text{sgn}{(#1)}} 
\newcommand{\SymPoly}[1]{\mathbf{m} (#1)} 
\newcommand{\UGHL}[1]{\mathbf{V}(#1)} 
\newcommand{\FGHL}[1]{\mathbf{L}(#1)}
\newcommand{\Ind}[3]{\text{Ind}_{#1}^{#2} \  #3}
\newcommand{\EditingNote}[1]{}
\newcommand{\PolyMod}[1]{\mathbf{W}(#1)}

\newcommand{\Partition}{{\mathbf P}}
\newcommand{\EpiArrow}{\twoheadrightarrow}
\newcommand{\SymFunc}[1]{\mathbf{A}_{#1}}

\newcommand{\Discriminant}[1]{\Omega_{#1}}
\newcommand{\Reg}[1]{#1^{\text{reg}}}
\newcommand{\C}[2]{{#1}(#2)}
\newcommand{\gc}{\mathrm{c}}
\newcommand{\gd}{\mathrm{d}}
\newcommand{\Singular}[1]{\mathbf{w}(#1)}
\newcommand{\sle}{\mathrm{e}}
\newcommand{\slf}{\mathrm{f}}
\newcommand{\slh}{\mathrm{h}}
\newcommand{\VermaGen}[1]{\mathrm{v}_{#1}}
\newcommand{\ie}{i.~e.~}
\newcommand{\cf}{cf.~}
\newcommand{\eg}{e.~g.~}
\newcommand{\Degree}{\text{deg} \hspace{0.1em}}
\newcommand{\ExpFunc}{\mathrm{Exp}}

\begin{document}
\begin{abstract}
Affine Lie algebras admit non-classical highest-weight theories through alternative partitions of the root system.
Although significant inroads have been made, much of the classical machinery is inapplicable in this broader context, and some fundamental questions remain unanswered.
In particular, the structure of the reducible objects in non-classical theories has not yet been fully understood.
This question is addressed here for affine $\SLTwo$, which has a unique non-classical highest-weight theory, termed ``imaginary''.
The reducible Verma modules in the imaginary theory possess an infinite descending series, with all factors isomorphic to a certain canonically associated module, the structure of which depends upon the highest weight.
If the highest weight is non-zero, then this factor module is irreducible, and conversely. 
This paper examines the degeneracy of the factor module of highest-weight zero.
The intricate structure of this module is understood via a realization in terms of the symmetric functions.
The realization permits the description of a family of singular (critical) vectors, and the classification of the irreducible subquotients.
The irreducible subquotients are characterized as those modules with an action given in terms of exponential functions, in the sense of Billig and Zhao.
\end{abstract}

\maketitle
\begin{section}{Overview}
Let $\K$ denote a field of characteristic zero and let $\g$ denote an affine Kac-Moody Lie algebra over $\K$.
Write $\CSA \subset \g$ for the Cartan subalgebra and $\RootSystem \subset \Dual{\CSA}$ for the root system.
Denote by $\g_{\phi}$ the root space associated to any root $\phi \in \RootSystem$.  So
$$ \g = (\oplus_{\phi \in \RootSystem} \g_\phi) \oplus \CSA, \qquad \Adjoint{h} |_{\g_\phi} = \phi (h), \quad h \in \CSA, \quad \phi \in \RootSystem.$$
A $\g$-module $V$ is called \NewTerm{weight} if the action of $\CSA$ upon $V$ is diagonalizable.  That is,
$$ V = \oplus_{\lambda \in \Dual{\CSA}} V_\lambda, \qquad h |_{V_\lambda} = \lambda (h), \quad h \in \CSA, \ \lambda \in \Dual\CSA.$$

\begin{subsection}{Highest-weight theories for affine Lie algebras}
The notion of a highest-weight module for $\g$ depends upon the partition of the root system.
A subset $\Partition \subset \RootSystem$ is called a \NewTerm{partition} of the root system if both
\begin{enumerate}
\item $\Partition$ is closed under root space addition, \ie 
if $\phi, \psi \in \Partition$ and $\phi + \psi \in \RootSystem$, then $\phi + \psi \in \Partition$;
\item $\Partition \cap -\Partition = \emptyset \ \text{and} \  \Partition \cup -\Partition = \RootSystem$.
\end{enumerate}
If $\RootSystem_+ (\pi)$ denotes the set of positive roots with respect to some basis $\pi \subset \RootSystem$ of the root system, then $\Partition = \RootSystem_+ (\pi)$ is an example of a partition.
A partition $\Partition$ defines a decomposition of the Lie algebra $\g$ as a direct sum of subalgebras
$$ \g = \NMinus \oplus \CSA \oplus \NPlus, \quad \text{where} \quad \NPlus = \oplus_{\phi \in \Partition}{\g_\phi}, \quad \NMinus = \oplus_{\phi \in \Partition}{\g_{-\phi}}.$$
A weight $\g$-module $V$ is of \NewTerm{highest-weight} $\lambda \in \Dual\CSA$ with respect to the partition $\Partition$ if there exists $v \in V_\lambda$ such that
$$ \UEA{\g} \cdot v = V, \quad \text{and} \quad \NPlus \cdot v = 0.$$
Thus the choice of partition $\Partition$ defines a theory of highest-weight modules.
A highest-weight theory defined by the set of all positive roots $\Partition = \RootSystem_{+} (\pi)$ with respect to some basis $\pi$ of the root system is called \NewTerm{classical}.
Let $\mathrm{W}$ denote the Weyl group associated to the root system $\RootSystem$.
Two partitions are \NewTerm{equivalent} if they are conjugate under the action of ${\mathrm W} \times \set{\pm 1}$.
Equivalent partitions define similar highest-weight theories.
All partitions of the root system of a finite-dimensional semisimple complex Lie algebra are classical, and hence equivalent.
In contrast, it has been shown by Jakobsen and Kac \cite{JK:ClassOfParts} and by Futorny \cite{F:ClassOfParts}, that there are finitely many, but never one, inequivalent partitions of the root system of an affine Lie algebra.
Highest-weight theories defined by inequivalent partitions are remarkably dissimilar (\cf subsection \ref{ImagHWTSLTwo}).
Thus any affine Lie algebra has multiple distinct highest-weight theories.
\end{subsection}

\begin{subsection}{Imaginary highest-weight theory for affine $\SLTwo$}\label{ImagHWTSLTwo}
There is precisely one non-classical partition, the \NewTerm{imaginary partition}, of the root system of affine $\SLTwo$.
The associated \NewTerm{imaginary highest-weight theory} has been pioneered by Futorny in \cite{F:ImagVerma}, \cite{F:Repns}.
These works provide an almost complete understanding of the universal objects of the theory, the \NewTerm{imaginary Verma modules}.

Let $\SLTwo = \SpanField{\K}{\sle,\slf,\slh}$ with the Lie bracket relations
$$ \Brac{\sle}{\slf} = \slh, \qquad \Brac{\slh}{\sle} = 2\sle, \qquad \Brac{\slh}{\slf} =-2\slf,$$
and Killing form $\langle \cdot,\cdot \rangle$, given by
\begin{eqnarray*}
\langle \slh,\slh \rangle = 2, \qquad \langle \sle,\slf \rangle = 1, \qquad
\langle \sle,\sle \rangle = \langle \slf,\slf \rangle = \langle \slh,\sle \rangle = \langle \slh,\slf \rangle = 0.
\end{eqnarray*}
Let $\g$ denote the affinization of $\SLTwo$:
\begin{equation}\label{BracketForG}
\g = {\SLTwo \otimes_{\K} \K[\mathrm{t}^{\pm1}]} \oplus {\K \gc} \oplus {\K \gd},
\end{equation}
with Lie bracket relations (writing $\C{x}{k}$ for ${x \otimes \mathrm{t}^k} \in \SLTwo \otimes \K [\mathrm{t}^{\pm1}]$):
\begin{eqnarray*}
\Brac{\C{x}{k}}{\C{y}{l}} &=& \C{\Brac{x}{y}}{k+l} + {k \delta_{k, -l} \langle x,y \rangle \gc}, \qquad \Brac{\gc}{\g} = 0, \\
\Brac{\gd}{\C{x}{k}} &=& k \C{x}{k}, \qquad x,y \in \SLTwo, \quad k,l \in \Z.
\end{eqnarray*}
The Cartan subalgebra $\CSA \subset \g$ is given by
$$ \CSA = \Span{\C{\slh}{0}, \gc, \gd}.$$
Let $\alpha, \delta \in \Dual{\CSA}$ be such that 
\begin{eqnarray*}
\alpha (\C{\slh}{0}) &=& 2, \qquad \alpha (\gc) = 0, \qquad \alpha(\gd) = 0, \\
\delta(\C{\slh}{0}) &=& 0, \qquad \delta(\gc) = 0, \qquad \delta (\gd) = 1.
\end{eqnarray*}
Then $\RootSystem = \set{\pm \alpha + i \delta | i \in \Z} \cup \set{i \delta | i \in \Z, \ i \ne 0}$.
The imaginary partition $\Partition \subset \RootSystem$ is given by
$$ \Partition = \set{ \alpha + i \delta | i \in \Z} \cup \set{ i \delta | i \in \Z, \ i > 0 }.$$
Thus the associated subalgebras $\NPlus, \NMinus$ of $\g$ are given by
$$ \NPlus = [ \oplus_{j \in \Z}{\K \C{\sle}{j}} ] \hspace{0.15em} \oplus \hspace{0.15em} [\oplus_{\substack{j > 0}}{\K \C{\slh}{j}}], \qquad \NMinus = [\oplus_{j \in \Z}{\K \C{\slf}{j}}] \hspace{0.15em} \oplus \hspace{0.15em} [\oplus_{\substack{j < 0}}{\K \C{\slh}{j}}].$$
All subsequent reference to a highest-weight module for $\g$ is with respect to the imaginary partition $\Partition$.
Let $\lambda \in \Dual \CSA$, and consider the one-dimensional vector space $\K \VermaGen{\lambda}$ as an $(\CSA \oplus \NPlus)$-module via
$$ \NPlus \cdot \VermaGen{\lambda} = 0, \qquad h \cdot \VermaGen{\lambda} = \lambda (h) \VermaGen{\lambda}, \quad h \in \CSA.$$
Let $\Verma{\lambda}$ denote the induced $\g$-module:
$$ \Verma{\lambda} = \Ind{\CSA \oplus \NPlus}{\g}{\K \VermaGen{\lambda}}.$$
The $\g$-module $\Verma{\lambda}$ is the universal highest-weight $\g$-module of highest-weight $\lambda$, and so is called an \NewTerm{imaginary Verma module}.
\begin{theoremnonum}\label{VermaStructureThm}\cite{F:ImagVerma}
Let $\lambda \in \Dual{\CSA}$.  Then:
\begin{enumerate}
\item If $\lambda (\gc) \ne 0$, then $\Verma{\lambda}$ is irreducible.
\item \label{CZeroCase} Suppose that $\lambda (\gc) = 0$.
Then $\Verma{\lambda}$ has an infinite descending series of submodules
\begin{equation*}
\Verma{\lambda} = V^0 \supset V^1 \supset V^2 \supset \cdots
\end{equation*}
such that any factor $V^i / V^{i+1}$, $i \geqslant 0$, is isomorphic to the quotient of $\g$-modules 
$$ \VermaQuotient{\lambda} = \Verma{\lambda} \  / \ { \langle \hspace{0.15em} \C{\slh}{j} \otimes \VermaGen{\lambda} \ | \  j < 0 \hspace{0.15em} \rangle },$$
up to a shift in the $\delta$ weight-decomposition.
Moreover, if $\lambda( \C{\slh}{0}) \ne 0$, then $\VermaQuotient{\lambda}$ is irreducible.
\end{enumerate}
\end{theoremnonum}
The value $\lambda (\gd)$ of the action of $\gd$ on the generating vector is immaterial to the structure of the imaginary Verma module $\Verma{\lambda}$.
Hence the Theorem \ref{VermaStructureThm} provides an almost complete description of the structure of the imaginary Verma modules for $\g$, lacking only a statement about the imaginary Verma module of highest-weight zero.
Part \ref{CZeroCase} of the Theorem motivates a study of $\Verma{0}$ through its canonically associated and infinitely occurrent quotient $\VermaQuotient{0}$.
This paper is an extensive study of the degeneracy of $\VermaQuotient{0}$.

Henceforth, $\g$ denotes the (unextended) loop algebra of $\SLTwo$,
\begin{equation*}
\g = {\SLTwo \otimes_{\K} \K[{\mathrm t}^{\pm1}]},
\end{equation*}
with Lie bracket relations
\begin{equation*}
\Brac{\C{x}{k}}{\C{y}{l}} = \C{\Brac{x}{y}}{k+l}, \qquad x,y \in \SLTwo, \quad k,l \in \Z.
\end{equation*}
and Cartan subalgebra $\CSA = \K \C{\slh}{0}$.
Degree in the indeterminate $\mathrm{t}$ endows $\g$ with a $\Z$-grading.
Let 
$$ \E = \oplus_{j \in \Z}{\K \C{\sle}{j}}, \qquad \Heis = \oplus_{j \in \Z}{\K \C{\slh}{j}}, \qquad \F = \oplus_{j \in \Z}{\K \C{\slf}{j}},$$
so that $\g = \E \oplus \Heis \oplus \F$ as a sum of graded abelian subalgebras.
All modules carry a compatible grading.
The $\g$-module $\VermaQuotient{0}$ may be defined by
\begin{equation*}
\VermaQuotient{0} = \Ind{\Heis \oplus \E}{\g}{\K \Gen}, \quad \text{where} \quad \Heis \oplus \E \cdot \Gen = 0.
\end{equation*}

For any weight $\g$-module $N$, write  
$$ N = \oplus_{\lambda \in \K}{N_\lambda}, \qquad \C{\slh}{0} |_{N_\lambda} = -2 \lambda, \quad \lambda \in \K,$$
for the decomposition of $N$ into eigenspaces for the operator $\C{\slh}{0}$.
An $\C{\slh}{0}$-eigenspace $N_\lambda$ of $N$ is called a \NewTerm{layer} of $N$, and the eigenspace decomposition is called the \NewTerm{layer decomposition} of $N$.
The layers of $N$ may be considered as modules for the subalgebra $\Heis$.
A \NewTerm{component} of $N$ is a graded component of some layer.
The $\g$-module $\VermaQuotient{0}$ has layer decomposition (\cf Corollary \ref{DecompositionCorollary})
$$ \VermaQuotient{0} = \oplus_{n \geqslant 0} \VermaQuotient{0}_n.$$
Write $\Level{n}$ for the layer $\VermaQuotient{0}_n$, $n \geqslant 0$.
\end{subsection}

\begin{subsection}{Layers and symmetric functions}
The $\Heis$-module structures of the layers $\Level{n}$ play a large part in the structural description of the $\g$-module $\VermaQuotient{0}$.
The $\Heis$-modules $\Level{n}$ have remarkable realizations in terms of the symmetric functions.
For any positive integer $n$, let 
$$\SymFunc{n} = \K[\mathrm{z}_1^{\pm1}, \dots, \mathrm{z}_{n}^{\pm1}]^{\SymGroup{n}}$$ 
denote the $\K$-algebra of symmetric Laurent polynomials in the $n$ indeterminates $\mathrm{z}_1, \dots, \mathrm{z}_n$.
It is shown in section \ref{SymFuncSection} that the layer $\Level{n}$ may be considered as a graded $\SymFunc{n}$-module in such a way that the action of $\Heis$ upon $\Level{n}$ factors through an epimorphism of algebras $ \UEA{\Heis} \EpiArrow \SymFunc{n}$ (\cf Theorem \ref{LevelSymFuncModuleThm}).
Therefore, it is sufficient to consider the layer $\Level{n}$ as a graded $\SymFunc{n}$-module.
Moreover, it is shown that as a graded $\SymFunc{n}$-module, the layer $\Level{n}$ is isomorphic to the graded regular module for $\SymFunc{n}$ (\cf Theorem \ref{LevelIsomorphismThm}).
Thus, in particular, the $\Heis$-module structure of $\Level{n}$ may be understood via the graded algebra $\SymFunc{n}$.
The classification of the irreducible quotients of the algebra $\SymFunc{n}$, and the construction of singular vectors of $\VermaQuotient{0}$, together permit the classification of the irreducible subquotients of the $\g$-module $\VermaQuotient{0}$.
\end{subsection}

\begin{subsection}{Singular vectors}
An element $v \in \VermaQuotient{0}$ is \NewTerm{singular} if $\E \cdot v = 0$.
If the layer $\Level{n}$ contains a non-zero singular vector, then the $\g$-module $\VermaQuotient{0}$ has submodule beginning on the layer $n$, and so the question of the existence of singular vectors is of principal interest.
A family of singular vectors is constructed on each layer $\Level{n}$ in section \ref{SingularVectorsSection}.
For any positive integer $n$, let
$$ \Discriminant{n} = \prod_{1 \leqslant i < j \leqslant n}{ (\mathrm{z}_i - \mathrm{z}_j)^2} \qquad \in \SymFunc{n}.$$
In the theory of symmetric functions, the function $\Discriminant{n}$ appears both as the square of the Vandermonde determinant, and as the discriminant function for degree-$n$ polynomials.
\begin{theoremnonum}
Let $n$ be a positive integer.  Then
\begin{enumerate}
\item All elements of the $\Heis$-submodule $\Discriminant{n} \cdot \Level{n}$ are singular.
\item The space $\Discriminant{n} \cdot \Level{n}$ is spanned by singular vectors of the form
$$
\Singular{\chi} = \sum_{\sigma \in \SymGroup{n}} \Sign{\sigma} \prod_{1 \leqslant i \leqslant n} \C{\slf}{\chi_i + \sigma(i)} \hspace{0.1em} \Gen
$$
where $\chi \in \Z^n$ and $\Gen$ denotes the generator of $\VermaQuotient{0}$.
\end{enumerate}
\end{theoremnonum}
\begin{conjecturenonum}
Let $n$ be a positive integer, and suppose that $v \in \Level{n}$ is singular.
Then $v \in \Discriminant{n} \cdot \Level{n}$.
\end{conjecturenonum}
The results of the present paper do not depend upon the truth of this conjecture.
\end{subsection}

\begin{subsection}{Irreducible quotients of the layers $\Level{n}$}
In section \ref{LevelISQSection}, the symmetric function realization is employed to classify the irreducible quotients of the $\Heis$-module $\Level{n}$.
Let $\varphi : \Z \to \K$ denote any function, and define
\begin{equation}\label{PolyModFunc}
\tilde{\varphi}: \UEA{\Heis} \to \K[ \mathrm{t}^{\pm 1}], \qquad \tilde{\varphi}: \C{\slh}{k} \mapsto \varphi(k) \mathrm{t}^k, \quad k \in \Z.
\end{equation}
Then $\tilde{\varphi}$ is a homomorphism of graded algebras.
Write $\PolyMod{\varphi}$ for $\Image{\tilde{\varphi}}$ considered as an $\Heis$-module via $\tilde{\varphi}$.
If $\K$ is algebraically closed and uncountable (\eg $\K = \FieldC$), then all irreducible $\Heis$-modules may be constructed in this way (\cf \cite{Chari}).
For any non-negative integer $n$, let $\ExpFunc_n$ denote the collection of all maps $\varphi: \Z \to \K$ such that
$$ \varphi (k) = -2 \sum_{i=1}^{n}{\alpha_i^k}, \quad k \in \Z,$$
for some $\alpha_1, \dots, \alpha_n \in \K$ all non-zero.  For example,
$$ \ExpFunc_0 = \set{0}, \qquad \ExpFunc_1 = \set{ \varphi : \Z \to \K | \varphi(k) = -2 \alpha^k \ \ \text{for some} \ \ \alpha \in \K, \ \alpha \ne 0}. $$
Let $\ExpFunc = \cup_{n \geqslant 0} \ExpFunc_n$.

\begin{theoremnonum}
For any non-negative integer $n$ and $\varphi \in \ExpFunc_n$, the $\Heis$-module $\PolyMod{\varphi}$ is an irreducible quotient of the $\Heis$-module $\Level{n}$. 
Moreover, if $\K$ is algebraically closed, then any irreducible quotient of the $\Heis$-module $\Level{n}$ is of the form $\PolyMod{\varphi}$ for some $\varphi \in \ExpFunc_n$.
\end{theoremnonum}
\end{subsection}

\begin{subsection}{Irreducible subquotients of the $\g$-module $\VermaQuotient{0}$}
A $\g$-module $Q$ is a \NewTerm{subquotient} of a $\g$-module $M$ if there exists a chain of $\g$-modules
$$ M \supset N \supset P $$
such that $N / P \cong Q$.
In section \ref{ISQSection}, the preceding results are employed to classify the irreducible subquotients of the $\g$-module $\VermaQuotient{0}$.
Associate to any weight $\Heis$-module $\Lambda$ an induced $\g$-module
$$ \UGHL{\Lambda} = \Ind{\Heis \oplus \E}{\g}{\Lambda} \qquad \text{where} \quad \E \cdot \Lambda = 0.$$
The $\g$-module $\UGHL{\Lambda}$ has a unique maximal submodule that has trivial intersection with $\Lambda$.
Denote by $\FGHL{\Lambda}$ the corresponding quotient of $\UGHL{\Lambda}$.
For any map $\varphi : \Z \to \K$, write $\FGHL{\varphi}$ for $\FGHL{\PolyMod{\varphi}}$.
\begin{theoremnonum}
For any $\varphi \in \ExpFunc$, the $\g$-module $\FGHL{\varphi}$ is an irreducible subquotient of $\VermaQuotient{0}$.
Moreover, if $\K$ is algebraically closed, then any irreducible subquotient of $\VermaQuotient{0}$ is of the form $\FGHL{\varphi}$ for some $\varphi \in \ExpFunc$.
\end{theoremnonum}

A function $\varphi : \Z \to \K$ is said to be \NewTerm{exp-polynomial} if there exists a non-negative integer $n$, polynomials $g_i \in \K[\mathrm{t}]$ and scalars $\alpha_i \in \K$, $1 \leqslant i \leqslant n$, such that
$$ \varphi(k) = \sum_{i=1}^{n}{g_i (k) \alpha_i^k}, \quad k \in \Z.$$
In \cite{BilligZhao}, it is shown that if $\varphi: \Z \to \K$ is exp-polynomial, then the components of the $\g$-module $\FGHL{\varphi}$ have finite dimension.
It follows that if $\K$ is algebraically closed, then the components of the irreducible subquotients of $\VermaQuotient{0}$ have finite dimension.
\end{subsection}

\end{section}

\begin{section}{The Canonical Quotient $\VermaQuotient{0}$}
The following preliminary result provides a description of the action of $\g$ upon $\VermaQuotient{0}$.
Here, and throughout, the use of a hat above a term indicates the omission of that term.
The subalgebra $\F$ is abelian, and so $\UEA{F}$ may be identified with the infinite-rank polynomial ring
$$ \K [ \hspace{0.1em} \C{\slf}{j} \hspace{0.1em} | \hspace{0.1em} j \in \Z \hspace{0.1em} ]. $$
\begin{propn}\label{FactsQuotientZero}
The following hold:
\begin{enumerate}
\item \label{FactsQuotientZero1} The $\g$-module $\VermaQuotient{0}$ is generated by an element $\Gen$ such that
the action of $\UEA{\F}$ on $\Gen$ is free, whilst the actions of $\UEA{\E}$ and $\UEA{\Heis}$ are trivial.
\item \label{UEAFAction} The $\g$-module $\VermaQuotient{0}$ has a basis
\begin{equation}\label{MonomialImages}
\cup_{n \geqslant 0} \set{ \FProd{n}{\gamma} \hspace{0.1em} \Gen | \gamma_1 \leqslant \gamma_2 \leqslant \cdots \leqslant \gamma_n, \ \gamma \in \Z^n}
\end{equation}
\item \label{AffineSLAction} The action of $\g$ on $\VermaQuotient{0}$ is given by
\begin{eqnarray*}
\C{\slf}{k} \cdot \FProd{n}{\gamma} \hspace{0.1em} \Gen &=&  \C{\slf}{k} \FProd{n}{\gamma} \hspace{0.1em} \Gen, \\
\C{\slh}{k} \cdot \FProd{n}{\gamma} \hspace{0.1em} \Gen &=& -2 \sum_{1 \leqslant i \leqslant n}{\C{\slf}{\gamma_1} \cdots \widehat{\C{\slf}{\gamma_i}} \cdots \C{\slf}{\gamma_n} \C{\slf}{\gamma_i + k} \hspace{0.1em} \Gen}, \\
\C{\sle}{k} \cdot \FProd{n}{\gamma} \hspace{0.1em} \Gen &=& -2 \sum_{1 \leqslant i < j \leqslant n}{\C{\slf}{\gamma_1} \cdots \widehat{\C{\slf}{\gamma_i}} \cdots \widehat{\C{\slf}{\gamma_j}} \cdots \C{\slf}{\gamma_n} \C{\slf}{\gamma_i + \gamma_j + k} \hspace{0.1em} \Gen}, 
\end{eqnarray*}
for all $\gamma \in \Z^n$, $n \geqslant 0$ and $k \in \Z$.
\end{enumerate}
\end{propn}
\begin{proof}
Part \ref{FactsQuotientZero1} is clear, and part \ref{UEAFAction} follows from part \ref{FactsQuotientZero1}.
To prove part \ref{AffineSLAction}, we firstly derive some commutation relations in $\UEA{\g}$ before considering them in light of parts \ref{FactsQuotientZero1} and \ref{UEAFAction}.
If $\mathcal{L}$ is any Lie algebra and $x \in \mathcal{L}$, then the adjoint map 
$$ \Adjoint{x} : \UEA{\mathcal{L}} \to \UEA{\mathcal{L}}, \quad \Adjoint{x} : y \mapsto \Brac{x}{y} = xy - yx, \quad y \in \UEA{\mathcal{L}}, $$
is a derivation of the associative product of $\UEA{\mathcal{L}}$.
That is,
\begin{equation}\label{AdjointDerivation}
\Brac{x}{\prod_{1 \leqslant i \leqslant n}{y_i}} = \sum_{1 \leqslant i \leqslant n}{y_1 \cdots y_{i-1} \Brac{x}{y_i} y_{i+1} \cdots y_{n}}.
\end{equation}
This formula yields immediately the commutation equation
\begin{equation*}
\Brac{\C{\slh}{k}}{\FProd{n}{\gamma}} = -2 \sum_{1 \leqslant i \leqslant n}{\C{\slf}{\gamma_1} \cdots \widehat{\C{\slf}{\gamma_i}} \cdots \C{\slf}{\gamma_n} \C{\slf}{\gamma_i + k} }   
\end{equation*}
for all $\gamma \in \Z^n$, $n \geqslant 0$ and $k \in \Z$.
Using formula \ref{AdjointDerivation} and substituting the above commutation equation for $\C{\slh}{k}$,  
\begin{eqnarray*}
\Brac{\C{\sle}{k}}{\FProd{n}{\gamma}}  &=& -2 \sum_{1 \leqslant i < j \leqslant n}{ \C{\slf}{\gamma_1} \cdots \widehat{\C{\slf}{\gamma_i}} \cdots \widehat{\C{\slf}{\gamma_j}} \cdots \C{\slf}{\gamma_r} \C{\slf}{\gamma_i + \gamma_j + k} } \\
	&& + \sum_{1 \leqslant i \leqslant n}{\C{\slf}{\gamma_1} \cdots \widehat{\C{\slf}{\gamma_i}} \cdots \C{\slf}{\gamma_r} \C{\slh}{\gamma_i + k} } 
\end{eqnarray*}
for all $\gamma \in \Z^n$, $n \geqslant 0$ and $k \in \Z$.
These formulae, in consideration of parts \ref{FactsQuotientZero1} and \ref{UEAFAction}, immediately imply the formulae of part \ref{AffineSLAction}, and completely describe the action of $\g$ on $\VermaQuotient{0}$.
\end{proof}

\begin{corollary}\label{DecompositionCorollary}
The $\g$-module $\VermaQuotient{0}$ has layer decomposition
$$ \VermaQuotient{0} = \oplus_{\substack{n \geqslant 0,\\ n \in \Z}}{\Level{n}},$$
where for any $n \geqslant 0$, 
$$ \Level{n} = \Span{\FProd{n}{\chi} \hspace{0.1em} \Gen | \chi \in \Z^n}.$$
\end{corollary}
\begin{remark}
It follows from Corollary \ref{DecompositionCorollary} above, and Proposition 3.6 of \cite{K:IDLA}, that the only integrable subquotient of $\VermaQuotient{0}$ is the trivial one-dimensional $\g$-module. 
\end{remark}
\end{section}

\begin{section}{Symmetric Function Realization of the Layers}\label{SymFuncSection}
This section presents a realization of the $\Heis$-module $\Level{n}$ as the graded regular module of the symmetric Laurent polynomials in $n$ variables.
The realization allows the classification of the irreducible quotients of the $\Heis$-module $\Level{n}$ outlined in section \ref{LevelISQSection}.

Fix a positive integer $n$.
The $\K$-algebra $\SymFunc{n}$ is $\Z$-graded by total degree.
The elementary symmetric functions $\ElemSymFunc{i} \in \SymFunc{n}$, $1 \leqslant i \leqslant n$, are defined by the polynomial equation
\begin{equation}\label{ElemSymFuncDefn}
\prod_{i=1}^n (1 + \mathrm{z}_i \mathrm{t}) = 1 + \sum_{i=1}^n \ElemSymFunc{i}(\mathrm{z}_1, \dots, \mathrm{z}_n) \hspace{0.1em} \mathrm{t}^i.
\end{equation}
Notice that $\ElemSymFunc{n} = \mathrm{z}_1 \cdots \mathrm{z}_n$ is invertible in $\SymFunc{n}$.
For any $k \in \Z$, let
$$ \PowerSum{k} = \mathrm{z}_1^k + \cdots + \mathrm{z}_n^k \quad \in \SymFunc{n},$$
denote the sum of $k$-powers of the indeterminates,
and for any $\gamma \in \Z^n$, write
$$
\SymPoly{\gamma} = \frac{1}{n!} \sum_{\sigma \in \SymGroup{n}}{\prod_{1 \leqslant i \leqslant n}{\mathrm{z}_{\sigma(i)}^{\gamma_i}}} \quad \in \SymFunc{n}.
$$
The symmetric polynomial $\SymPoly{\gamma}$ may alternatively be defined by
$$
\SymPoly{\gamma} = \frac{1}{n!} \sum_{\sigma \in \SymGroup{n}}{\prod_{1 \leqslant i \leqslant n}{ \mathrm{z}_i^{\gamma_{\sigma(i)}} } }.
$$
The set
$\set{ \SymPoly{\gamma} | \gamma \in \Z^n }$
spans the $\K$-algebra $\SymFunc{n}$ of symmetric Laurent polynomials.

\begin{lemma}\label{MultOfSymPolys}
Let $n > 0 $ and $\gamma, \chi \in \Z^n$.  Then
$$
\SymPoly{\gamma} \cdot \SymPoly{\chi} = \frac{1}{n!} \sum_{\tau \in \SymGroup{n}}{\SymPoly{\gamma + \chi_\tau}},
$$
where $\chi_\tau \in \Z^n$ is given by $(\chi_\tau)_i = \chi_{\tau (i)}$, for all $1 \leqslant i \leqslant n$ and $\tau \in \SymGroup{n}$.
\end{lemma}
\begin{proof}
Let $\gamma, \chi \in \Z^n$.  Then
\begin{eqnarray*}
\SymPoly{\gamma} \cdot \SymPoly{\chi} &=& \big(\frac{1}{n!}\big)^2 \sum_{\sigma \in \SymGroup{n}} \sum_{\tau \in \SymGroup{n}} \prod_{1 \leqslant i \leqslant n} \mathrm{z}_i^{\gamma_{\sigma (i)} + \chi_{\tau (i)}} \\
 &=& \big(\frac{1}{n!}\big)^2 \sum_{\sigma \in \SymGroup{n}} \sum_{\tau \in \SymGroup{n}} \prod_{1 \leqslant i \leqslant n} \mathrm{z}_i^{\gamma_{\sigma (i)} + \chi_{(\tau \circ \sigma)(i)}} \\
 && \Note{substituting $\tau \circ \sigma$ for $\tau$} \\
 &=& \big(\frac{1}{n!}\big)^2 \sum_{\sigma \in \SymGroup{n}} \sum_{\tau \in \SymGroup{n}} \prod_{1 \leqslant i \leqslant n} \mathrm{z}_i^{(\gamma + \chi_\tau)_{\sigma (i)}} \\
 && \Note{since $\chi_{\tau \circ \sigma} = (\chi_\tau)_\sigma$} \\
 &=& \frac{1}{n!} \sum_{\tau \in \SymGroup{n}} \SymPoly{\gamma + \chi_\tau} 
\end{eqnarray*}
\end{proof}
\begin{propn}\label{PowerSumsGenerate}
For any positive integer $n$:
\begin{enumerate}
\item\label{PSGPart1} The $\K$-algebra $\SymFunc{n}$ is generated by the set $\set{\ElemSymFunc{i} | 1 \leqslant i \leqslant n} \cup \set{ \Inverse{\ElemSymFunc{n}}}$.
\item\label{PSGPart2} The $\K$-algebra $\SymFunc{n}$ is generated by the set of power sums $\set{\PowerSum{k} | k \in \Z}$.
\end{enumerate}
\end{propn}
\begin{proof}
Let 
$$\SymFunc{n}^{+} = \K [ \mathrm{z}_1, \dots, \mathrm{z}_n ]^{\SymGroup{n}}, \qquad \SymFunc{n}^{-} = \K [ \mathrm{z}_1^{-1}, \dots, \mathrm{z}_n^{-1} ]^{\SymGroup{n}}.$$
Then
\begin{equation}\label{PSGEq1}
\SymFunc{n} = \sum_{k \ge 0}{\ElemSymFunc{n}^{-k} \cdot \SymFunc{n}^{+}},
\end{equation}
and, in particular,
\begin{equation}\label{PSQEq2}
\SymFunc{n} = \SymFunc{n}^{-} \cdot \SymFunc{n}^{+}.
\end{equation}
Part \ref{PSGPart1} follows from the Fundamental Theorem of Symmetric Functions and equation \ref{PSGEq1},
while part \ref{PSGPart2} follows from the Newton-Girard formulae and equation \ref{PSQEq2}.
\end{proof}
\begin{propn}
For any $n > 0$, the layer $\Level{n}$ is a $\Z$-graded $\SymFunc{n}$-module via linear extension of
\begin{equation}\label{SymFuncAction}
\SymPoly{\gamma} \cdot \FProd{n}{\chi} \hspace{0.1em} \Gen = \frac{1}{n!} \sum_{\sigma \in \SymGroup{n}} \prod_{1 \leqslant i \leqslant n}{ \C{\slf}{\chi_i + \gamma_{\sigma (i)}}} \hspace{0.1em} \Gen, \qquad \gamma, \chi \in \Z^n,
\end{equation}
with the $\Z$-grading defined by
$$ \Degree {\FProd{n}{\chi} \hspace{0.1em} \Gen } = \sum_{1 \leqslant i \leqslant n}{\chi_i}, \qquad \chi \in \Z^n.$$
\end{propn}
\begin{proof}
For any $\varepsilon, \gamma, \chi \in \Z^n$,
\begin{eqnarray*}
\SymPoly{\varepsilon} \cdot \big( \SymPoly{\gamma} \cdot \FProd{n}{\chi} \hspace{0.1em} \Gen \big) 
	&=& \big(\frac{1}{n!}\big)^2 \sum_{\sigma \in \SymGroup{n}} \sum_{\tau \in \SymGroup{n}} \prod_{1 \leqslant i \leqslant n}{ \C{\slf}{\chi_i + \varepsilon_{\sigma (i)} + \gamma_{\tau (i)}}} \hspace{0.1em} \Gen \\
	&=& \frac{1}{n!} \sum_{\tau \in \SymGroup{n}} \frac{1}{n!} \sum_{\sigma \in \SymGroup{n}} \prod_{1 \leqslant i \leqslant n}{ \C{\slf}{\chi_i + (\varepsilon + \gamma_\tau)_{\sigma (i)}}} \hspace{0.1em} \Gen \\
	&& \Note{substituting $\tau \circ \sigma$ for $\tau$} \\
	&=& \Big( \frac{1}{n!} \sum_{\tau \in \SymGroup{n}} \SymPoly{\varepsilon + \gamma_\tau} \Big) \cdot \FProd{n}{\chi} \hspace{0.1em} \Gen \\
	&=& \big( \SymPoly{\varepsilon} \SymPoly{\gamma} \big) \cdot \FProd{n}{\chi} \hspace{0.1em} \Gen,
\end{eqnarray*}
by Lemma \ref{MultOfSymPolys}.
As the polynomials $\SymPoly{\gamma}$ span $\SymFunc{n}$, linear extension of \ref{SymFuncAction} endows $\Level{n}$ with the structure of a $\Z$-graded $\SymFunc{n}$-module.
\end{proof}
\begin{theorem}\label{LevelSymFuncModuleThm}
Let $n > 0$.
The action of $\UEA{\Heis}$ on $\Level{n}$ factors through an epimorphism
$$ \Psi: \UEA{\Heis} \EpiArrow \SymFunc{n}$$
of graded algebras defined by
$$ \Psi: \C{\slh}{k} \mapsto -2 \PowerSum{k}, \qquad k \in \Z.$$
That is, if $\rho$ and $\nu$ denote the representations of $\UEA{\Heis}$ and $\SymFunc{n}$ on $\Level{n}$, respectively, then the following diagram commutes:
$$
\xymatrix{
\UEA{\Heis} \ar@{->>}[r]^-\Psi \ar[rd]_\rho& \SymFunc{n}\ar[d]^\nu \\
&\End{\Level{n}}
}
$$
\end{theorem}
\begin{proof}
The map $\Psi$ is an algebra epimorphism by Proposition \ref{PowerSumsGenerate} part \ref{PSGPart2}.
Let $k \in \Z$, and let
$ \iota_k = (k,0, \dots,0) \in \Z^n$.
Then
$$ \PowerSum{k} = n \SymPoly{\iota_k}.$$
It follows that, for any $\chi \in \Z^n$,
$$ \PowerSum{k} \cdot \FProd{n}{\chi} \hspace{0.1em} \Gen 
	= \sum_{1 \leqslant i \leqslant n}{ \C{\slf}{\chi_1} \cdots \widehat{\C{\slf}{\chi_i}} \cdots \C{\slf}{\chi_n} \C{\slf}{\chi_i + k} } \hspace{0.1em} \Gen $$
and so
$ \C{\slh}{k} |_{\Level{n}} = \Psi( \C{\slh}{k}) |_{\Level{n}} $
by Proposition \ref{FactsQuotientZero} part \ref{AffineSLAction}.
\end{proof}
Therefore, for $n > 0$, it is sufficient to consider $\Level{n}$ as a $\Z$-graded $\SymFunc{n}$-module.
Write $\Reg{\SymFunc{n}}$ for the regular $\Z$-graded $\SymFunc{n}$-module, \ie for $\SymFunc{n}$ considered as a $\Z$-graded $\SymFunc{n}$-module under multiplication.
\begin{theorem}\label{LevelIsomorphismThm}
For any $n > 0$, the map
$$ \Theta : \Level{n} \to \Reg{\SymFunc{n}} $$
defined by linear extension of
$$ \Theta: \FProd{n}{\chi} \hspace{0.1em} \Gen \mapsto \SymPoly{\chi}, \qquad \chi \in \Z^n,$$
is an isomorphism of $\Z$-graded $\SymFunc{n}$-modules.
\end{theorem}
\begin{proof}
The map $\Theta$ is a bijection, by Proposition \ref{FactsQuotientZero} part \ref{UEAFAction}.
Let $\gamma, \chi \in \Z^n$.  Then
\begin{eqnarray*}
\Theta \Big( \SymPoly{\gamma} \cdot \FProd{n}{\chi} \hspace{0.1em} \Gen \Big) 
	&=& \Theta \Big( \frac{1}{n!} \sum_{\sigma \in \SymGroup{n}} \prod_{1 \leqslant i \leqslant n}{ \C{\slf}{\chi_i + \gamma_{\sigma (i)}}} \hspace{0.1em} \Gen \Big) \\
	&=& \frac{1}{n!} \sum_{\sigma \in \SymGroup{n}} \SymPoly{\chi + \gamma_\sigma} \\
	&=& \SymPoly{\gamma} \cdot \SymPoly{\chi} \ENote{by Lemma \ref{MultOfSymPolys}} \\
	&=& \SymPoly{\gamma} \cdot \Theta \Big( \FProd{n}{\chi} \hspace{0.1em} \Gen \Big).
\end{eqnarray*}
Hence $\Theta$ is an isomorphism of $\Z$-graded $\SymFunc{n}$-modules.
\end{proof}
\begin{corollary}\label{UHvIsoMr}
Let $n \geqslant 0$ and let $v \in \Level{n}$ be non-zero and homogeneous.
Then $\UEA{\Heis}v$ and $\Level{n}$ are isomorphic as $\Heis$-modules.
\end{corollary}
\begin{proof}
For $n=0$, the statement is trivial.
For $n > 0$, it is sufficient to employ Theorem \ref{LevelIsomorphismThm}, and observe that the ring $\SymFunc{n}$ is an integral domain.
\end{proof}
\end{section}

\begin{section}{Singular Vectors}\label{SingularVectorsSection}
The existence of non-zero singular vectors is related to the degeneracy of the module $\VermaQuotient{0}$: 
\begin{propn}\label{SVPropn}
Let $n \geqslant 0$.  Then
\begin{enumerate}
\item \label{SVPropn1} The set of all singular vectors in $\Level{n}$ form an $\Heis$-submodule.
\item \label{SVPropn2} If $v \in \Level{n}$ is non-zero and singular, then $\VermaQuotient{0}$ has a $\g$-submodule beginning on layer $n$.
That is, if $V = \UEA{\g}v$, then $V$ has layer decomposition
$$ V = \oplus_{m \geqslant n}{V_m}, \qquad V_n \ne 0.$$
\end{enumerate}
\end{propn}
\begin{proof}
The set of all singular vectors in $\Level{n}$ clearly form a vector space.
Now suppose that $v \in \Level{n}$ is singular.  Then
$$ \C{\sle}{k} (\C{\slh}{l} v) = 2 \C{\sle}{k+l} v + \C{\slh}{l} \C{\sle}{k} v = 0,$$
for any $k,l \in \Z$.
Thus if $v$ is singular, then so is $\C{\slh}{l} v$, for any $l \in \Z$, and so the set of singular vectors in $\Level{n}$ form an $\Heis$-module, proving part \ref{SVPropn1}.

For part \ref{SVPropn2}, suppose again that $v \in \Level{n}$ is a non-zero singular vector,
and let $V = \UEA{\g}v$.  Then
$$ V = \UEA{\g} \cdot v = \UEA{\F} \otimes \UEA{\Heis} \cdot v \subset \oplus_{m \geqslant n}{\Level{m}},$$
since
$\UEA{\g} = \UEA{\F} \otimes \UEA{\Heis} \otimes \UEA{\E}$. 
\end{proof}

\begin{theorem}
For any $n > 0$ and $\chi \in \Z^n$,
$$\Singular{\chi} = \sum_{\sigma \in \SymGroup{n}}{ \Sign{\sigma} \prod_{1 \leqslant i \leqslant n}{ \C{\slf}{\chi_i + \sigma(i)}} } \hspace{0.1em} \Gen \qquad \in \Level{n}$$
is a singular vector.
\end{theorem}
\begin{proof}
Singularity may be demonstrated directly by applying the formula for the action of $\C{\sle}{k}$, $k \in \Z$, of Proposition \ref{FactsQuotientZero} part \ref{AffineSLAction}.
\end{proof}

\begin{lemma}\label{DiscriminantLemma}
For any $n > 0$, the symmetric function $\Discriminant{n}$ is equal to
$$ \sum_{\sigma, \tau \in \SymGroup{n}} \Sign{\sigma \circ \tau} \prod_{1 \leqslant i \leqslant n} \mathrm{z}_i^{\sigma(i) + \tau(i) - 2}$$
up to a change in sign.
\end{lemma}
\begin{proof}
Let $\Phi_n = \sum_{\sigma \in \SymGroup{n}} \Sign{\sigma} \prod_{1 \leqslant i \leqslant n} \mathrm{z}_i^{\sigma (i)}$.
It is not difficult to demonstrate that 
$$\Phi_n |_{\mathrm{z}_i = \mathrm{z}_j} = 0, \qquad 1 \leqslant i < j \leqslant n,$$ 
and that $\Phi_n |_{\mathrm{z}_i} = 0$ for all $1 \leqslant i \leqslant n$.
Hence $\Phi_n$ is equal up to sign to
$$ \prod_{1 \leqslant i \leqslant n} \mathrm{z}_i\ \cdot \ \prod_{1 \leqslant i < j \leqslant n} (\mathrm{z}_i - \mathrm{z}_j),$$
by degree considerations.
Therefore $\Discriminant{n}$ and
$$ \prod_{1 \leqslant i \leqslant n} \mathrm{z}_i^{-2} \ \cdot \ \Phi_n^2 $$
are equal up to sign, from which the result follows.
\end{proof}

\begin{theorem}\label{ExistenceSV}
For any $n > 0$,
$$ \Discriminant{n} \cdot \Level{n} = \Span{\Singular{\chi} | \chi \in \Z^n},$$
and hence all elements of the non-zero $\Heis$-submodule $\Discriminant{n} \cdot \Level{n}$ are singular.
\end{theorem}
\begin{proof}
Fix $n > 0$, and let
$$ W_n = \Span{\Singular{\chi} | \chi \in \Z^n}.$$
The symmetric function realization of Theorem \ref{LevelIsomorphismThm} may be used to demonstrate the inclusion $W_n \subset \Discriminant{n} \cdot \Level{n}$.
Let $\chi \in \Z^n$.  Then
\begin{eqnarray*}
\Theta (\Singular{\chi}) &=& \frac{1}{n!} \sum_{\tau \in \SymGroup{n}} \sum_{\sigma \in \SymGroup{n}} \Sign{\sigma} \prod_{1 \leqslant i \leqslant n} \mathrm{z}_i^{\chi_{\tau (i)} + \sigma(\tau (i))} \\
	&=& \frac{1}{n!} \sum_{\tau \in \SymGroup{n}} \Sign{\tau} \sum_{\sigma \in \SymGroup{n}} \Sign{\sigma} \prod_{1 \leqslant i \leqslant n} \mathrm{z}_i^{\chi_{\tau(i)} + \sigma(i)} \\
	&& \Note{substituting $\sigma \circ \tau^{-1}$ for $\sigma$} \\
	&=& \frac{1}{n!} \sum_{\tau \in \SymGroup{n}} \Sign{\tau} F_\tau,
\end{eqnarray*}
where, for any $\tau \in \SymGroup{n}$,
$$ F_\tau = \sum_{\sigma \in \SymGroup{n}} \Sign{\sigma} \prod_{1 \leqslant i \leqslant n} \mathrm{z}_i^{\chi_{\tau(i)} + \sigma(i)}.$$
It is not difficult to verify that 
$$ F_\tau |_{\mathrm{z}_i = \mathrm{z}_j} = 0, \qquad 1 \leqslant i < j \leqslant n, \quad \tau \in \SymGroup{n}.$$
Therefore $\Theta(\Singular{\chi})$ is divisible by $\Discriminant{n}$ in $\SymFunc{n}$.  That is,
$ \Theta(\Singular{\chi}) \in \Discriminant{n} \cdot \Level{n},$
and hence $W_n \subset \Discriminant{n} \cdot \Level{n}$.

Now let $\Lambda_n = (\prod_{1 \leqslant i \leqslant n} \mathrm{z}_i^2) \cdot \Discriminant{n}.$
The factor $\prod_{1 \leqslant i \leqslant n}{\mathrm{z}_i^2}$ is invertible in $\SymFunc{n}$, and so 
$$ \Discriminant{n} \cdot \Level{n} = \Lambda_n \cdot \Level{n},$$
by Theorem \ref{LevelIsomorphismThm}.
In particular, by Corollary \ref{DecompositionCorollary},
$$ \Discriminant{n} \cdot \Level{n} = \Span{\Lambda_n \cdot \FProd{n}{\chi} \hspace{0.1em} \Gen | \chi \in \Z^n}.$$
By Lemma \ref{DiscriminantLemma}, the polynomial $\Lambda_n$ is equal up to sign to 
$$ \sum_{\sigma, \tau \in \SymGroup{n}} \Sign{\sigma \circ \tau} \prod_{1 \leqslant i \leqslant n} \mathrm{z}_i^{\sigma(i) + \tau(i)}.$$
Therefore, for any $\chi \in \Z^n$,
\begin{eqnarray*}
\Lambda_n \cdot \FProd{n}{\chi} \hspace{0.1em} \Gen
	&=& \frac{1}{n!} \sum_{\nu \in \SymGroup{n}} \sum_{\sigma, \tau \in \SymGroup{n}} \Sign{\sigma \circ \tau} \prod_{1 \leqslant i \leqslant n} \C{\slf}{\chi_i + (\sigma \circ \nu)(i) + (\tau \circ \nu)(i)} \hspace{0.1em} \Gen \\
	&=& \frac{1}{n!} \sum_{\nu \in \SymGroup{n}} \sum_{\sigma, \tau \in \SymGroup{n}} \Sign{\sigma \circ \tau} \prod_{1 \leqslant i \leqslant n} \C{\slf}{\chi_i + \sigma(i) + \tau(i)} \hspace{0.1em} \Gen \\
	&& \Note{substituting $\sigma \circ \nu^{-1}$ for $\sigma$ and $\tau \circ \nu^{-1}$ for $\tau$} \\
	&=& \sum_{\sigma \in \SymGroup{n}} \Sign{\sigma} \sum_{\tau \in \SymGroup{n}} \Sign{\tau} \prod_{1 \leqslant i \leqslant n} \C{\slf}{\chi_i + \sigma(i) + \tau(i)} \hspace{0.1em} \Gen \\
	&=& \sum_{\sigma \in \SymGroup{n}} \Sign{\sigma} \hspace{0.1em} \Singular{\chi(\sigma)},
\end{eqnarray*}
where $\chi(\sigma) \in \Z^n$ is given by
$$ \chi(\sigma)_i = \chi_i + \sigma(i), \qquad 1 \leqslant i \leqslant n,$$
for all $\sigma \in \SymGroup{n}$.
Hence $\Discriminant{n} \cdot \Level{n} \subset W_n$.
\end{proof}
\begin{conjecture}
Let $n >0$, and suppose that $v \in \Level{n}$ is singular.
Then $v \in \Discriminant{n} \cdot \Level{n}$.
\end{conjecture}
\end{section}

\begin{section}{Irreducible Quotients of the Layers}\label{LevelISQSection}
For any $\Z$-graded $\K$-algebra $B$, write $B^{(k)}$, $k \in \Z$, for the graded components of $B$.
\begin{propn}\label{ImageShapePropn}
Let $n$ be a positive integer, and let $B$ be a simple quotient of the graded algebra $\SymFunc{n}$.
Then $B = \FieldF [ \mathrm{t}^{\pm m}]$ for some positive divisor $m$ of $n$ and finite algebraic field extension $\FieldF$ of $\K$.
\end{propn}
\begin{proof}
As $\ElemSymFunc{n} = \mathrm{z}_1 \cdots \mathrm{z}_n$ is invertible in $\SymFunc{n}$, it must be that $B^{(n)} \ne 0$.
Let $m$ be the minimal positive integer such that $B^{(m)} \ne 0$, and let $u \in B^{(m)}$ be non-zero.
Then $u$ is invertible, since $B$ is simple, and so multiplication by $u^k$ is a vector-space automorphism of $B$ such that
$$ B^{(l)} \to B^{(km + l)}, \qquad l \in \Z. $$
In particular,
\begin{equation}\label{AutoDecomp}
B^{(km)} = u^k \cdot B^{(0)}, \qquad k \in \Z.
\end{equation}
Suppose that $B^{(l)} \ne 0$, for some $l \in \Z$, and let $q,r$ be the unique integers such that
$$ l = qm + r, \qquad 0 \leqslant r < m.$$
Then
$$ 0 \ne u^{-qm} \cdot B^{(l)} = B^{(r)},$$
and so $r=0$ by the minimality of $m$.  Hence
$$ B = \oplus_{k \in \Z}{ B^{(km)} },$$
and in particular $m$ must be a divisor of $n$.
Moreover, by \ref{AutoDecomp},
$$ B \cong B^{(0)} \otimes_{\K} \K[ \mathrm{t}^{\pm 1} ]$$
via $u^k \mapsto \mathrm{t}^{km}$, $k \in \Z$.
As $\SymFunc{n}$ is finitely generated, by Proposition \ref{PowerSumsGenerate} part \ref{PSGPart1}, so is the $\K$-algebra $\SymFunc{n}^{(0)}$.
Hence $B^{(0)}$ is a finite algebraic field extension of $\K$ (see, for example, \cite{AtiyahMacdonald}, Proposition 7.9).
\end{proof}

\begin{propn}\label{HomClassificationPropn}
Let $n$ be a positive integer, and suppose that $\zeta : \SymFunc{n} \to \K$ is a non-zero algebra homomorphism.
Then there exist scalars $\alpha_1, \dots, \alpha_n$, all non-zero and algebraic over $\K$, such that
$$ \zeta (\PowerSum{k}) = \sum_{i=1}^n \alpha_i^k, \qquad k \in \Z.$$
\end{propn}
\begin{proof}
Suppose that $\zeta: \SymFunc{n} \to \K$ is a non-zero homomorphism, and let 
$$ g(\mathrm{t}) = 1 + \sum_{i=1}^n \zeta(\ElemSymFunc{i}) \hspace{0.1em} \mathrm{t}^i \quad \in \K [ \mathrm{t} ].$$
As $\ElemSymFunc{n} = \mathrm{z}_1 \cdots \mathrm{z}_n$ is invertible in $\SymFunc{n}$, it must be that $\zeta(\ElemSymFunc{n}) \ne 0$.
Let $\alpha_1, \dots \alpha_n$ be some iteration of the scalars defined by
$$ g (\mathrm{t}) = \prod_{i=1}^n (1 + \alpha_i \mathrm{t}).$$
Then by equation \ref{ElemSymFuncDefn},
\begin{equation}\label{EvalProp}
\zeta (\ElemSymFunc{i}) = \ElemSymFunc{i}(\alpha_1, \dots, \alpha_n), \qquad 1 \leqslant i \leqslant n.
\end{equation}
The $\alpha_i$ are necessarily non-zero since
$$ \alpha_1 \cdots \alpha_n = \ElemSymFunc{n}(\alpha_1, \dots, \alpha_n) = \zeta(\ElemSymFunc{n}) \ne 0.$$
By Proposition \ref{PowerSumsGenerate} part \ref{PSGPart1}, there is a unique algebra homomorphism with the property \ref{EvalProp}, namely the restriction of the evaluation map
$$ \mathrm{z}_i \mapsto \alpha_i, \qquad 1 \leqslant i \leqslant n.$$
In particular, $\zeta (\PowerSum{k}) = \sum_{i=1}^n \alpha_i^k$, for all $k \in \Z$.
\end{proof}

\begin{theorem}\label{LevelImageClassification}
For any non-negative integer $n$ and $\varphi \in \ExpFunc_n$, the $\Heis$-module $\PolyMod{\varphi}$ is an irreducible quotient of the $\Heis$-module $\Level{n}$. 
Moreover, if $\K$ is algebraically closed, then any irreducible quotient of the $\Heis$-module $\Level{n}$ is of the form $\PolyMod{\varphi}$ for some $\varphi \in \ExpFunc_n$.
\end{theorem}
\begin{proof}
The statement is trivial for $n=0$, so suppose that $n$ is a positive integer.
Let $\varphi \in \ExpFunc_n$, and let $\alpha_1, \dots \alpha_n \in \K$ be non-zero scalars such that
$$ \varphi (k) = -2 \sum_{i=1}^{n}{\alpha_i^k}, \qquad k \in \Z.$$
Define 
$$\eta = \eta_{\alpha_1, \dots, \alpha_n}: \K [ \mathrm{z}_1^{\pm 1}, \dots, \mathrm{z}_n^{\pm 1}] \to \K [ \mathrm{t}^{\pm 1}]$$
by extension of $\mathrm{z}_i \mapsto \alpha_i \mathrm{t}$, $1 \leqslant i \leqslant n$.
Then
$$ \eta |_{\SymFunc{n}} : \SymFunc{n} \to \Image{\eta |_{\SymFunc{n}}} \subset \K [ \mathrm{t}^{\pm 1}]$$
is a homomorphism of graded algebras, and so $\Image{\eta |_{\SymFunc{n}}}$ may be considered as an $\SymFunc{n}$-module, and as a quotient  of the regular $\SymFunc{n}$-module $\Reg{\SymFunc{n}}$.
Therefore $\Image{\eta |_{\SymFunc{n}}}$ is a quotient of the $\SymFunc{n}$-module $\Level{n}$, by Theorem \ref{LevelIsomorphismThm}.
Now $\Image{\eta |_{\SymFunc{n}}}$ is an $\Heis$-module via the map $ \Psi: \UEA{\Heis} \to \SymFunc{n}$ of Theorem \ref{LevelSymFuncModuleThm}.
For any $k \in \Z$,
\begin{eqnarray*}
( \eta |_{\SymFunc{n}} \circ \Psi ) ( \C{\slh}{k} ) &=& -2 \eta |_{\SymFunc{n}} (\PowerSum{k}) \\
	&=& -2 \sum_{i=1}^n \alpha_i^k \hspace{0.1em} \mathrm{t}^k \\
	&=& \tilde{\varphi} (\C{\slh}{k}),
\end{eqnarray*}
where $\tilde{\varphi} : \UEA{H} \to \K [ \mathrm{t}^{\pm 1}]$ is defined by \ref{PolyModFunc}.
Hence $\Image{\eta |_{\SymFunc{n}}} \cong \PolyMod{\varphi}$ as $\Heis$-modules, and so $\PolyMod{\varphi}$ is a quotient of the $\Heis$-module $\Level{n}$.
It is not too difficult to verify that $\PolyMod{\varphi}$ is an irreducible $\Heis$-module for any $\varphi \in \ExpFunc$.

Now suppose that $\K$ is algebraically closed, and that $\Gamma$ is an irreducible quotient of the $\Heis$-module $\Level{n}$.
By Theorems \ref{LevelSymFuncModuleThm} and \ref{LevelIsomorphismThm}, $\Gamma$ is an $\SymFunc{n}$-module, and a quotient of $\Reg{\SymFunc{n}}$.
Hence there exists a simple quotient $B$ of the graded algebra $\SymFunc{n}$
$$ \eta : \SymFunc{n} \EpiArrow B,$$
such that $B \cong \Gamma$ as $\SymFunc{n}$-modules, when $B$ is considered as an $\SymFunc{n}$-module via the algebra epimorphism $\eta$.
Moreover, by Proposition \ref{ImageShapePropn}, $B = \K [\mathrm{t}^{\pm m}]$ for some positive divisor $m$ of $n$.
Let $\zeta : \SymFunc{n} \to \K$ be given by 
$$ \eta ( x ) = \zeta (x) \hspace{0.1em} \mathrm{t}^{\Degree{x}},$$
for all homogeneous $x \in \SymFunc{n}$.
Then $\zeta$ is a non-zero algebra homomorphism, and so by Proposition \ref{HomClassificationPropn} there exist non-zero scalars $\alpha_1, \dots, \alpha_n \in \K$ such that 
$$ \eta (\PowerSum{k}) = \zeta (\PowerSum{k}) \hspace{0.1em} \mathrm{t}^k = \sum_{i=1}^n \alpha_i^k \hspace{0.1em} \mathrm{t}^k, \qquad k \in \Z.$$
Therefore $\Gamma \cong \PolyMod{\varphi}$, where
$$\varphi : \Z \to \K, \qquad \varphi (k) = -2 \sum_{i=1}^n \alpha_i^k, \qquad k \in \Z,$$
by Theorem \ref{LevelSymFuncModuleThm}.
\end{proof}
\end{section}

\begin{section}{Irreducible Subquotients of $\VermaQuotient{0}$}\label{ISQSection}
A layer $\Lambda$ of a weight $\g$-module $N$ is said to be a \NewTerm{generative highest layer} if
$$ \UEA{\g} \cdot \Lambda = N, \quad \text{and} \quad \E \cdot \Lambda = 0.$$
Any two generative highest layers of $N$ must coincide, and so we may speak of the generative highest layer, if one exists.
For any weight $\Heis$-module $\Lambda$, write $\UGHL{\Lambda}$ for the induced $\g$-module with generative highest layer $\Lambda$, i.~e.
$$ \UGHL{\Lambda} = \Ind{\Heis + \E}{\g}{\Lambda}, \quad \text{where} \quad \E \cdot \Lambda = 0.$$
\begin{propn}
Let $\Lambda$ be a weight $\Heis$-module.  Then the $\g$-module $\UGHL{\Lambda}$ has a unique maximal submodule that has trivial intersection with $\Lambda$.
\end{propn}
\begin{proof}
Let $\lambda \in \K$ be such that $\C{\slh}{0} |_\Lambda = -2 \lambda$.  
Then $\UGHL{\Lambda}$ has layer decomposition
$$ \UGHL{\Lambda} = \oplus_{n \geqslant 0}{\UGHL{\Lambda}_{\lambda + n}},$$
and $\UGHL{\Lambda}_\lambda = \Lambda$.
If $N \subset \UGHL{\Lambda}$ is a $\g$-submodule that has trivial intersection with $\Lambda$, then
$$ N \subset \oplus_{n > 0}{\UGHL{\Lambda}_{\lambda + n}}.$$
Hence the same is true of the sum of all such $\g$-submodules.
This sum is itself a $\g$-submodule, and its maximality and uniqueness follow from construction.
\end{proof}
For any weight $\Heis$-module $\Lambda$, denote by $\FGHL{\Lambda}$ the quotient of $\UGHL{\Lambda}$ by its unique maximal submodule that has trivial intersection with $\Lambda$.
Hence $\FGHL{\Lambda}$ is an irreducible $\g$-module if $\Lambda$ is an irreducible $\Heis$-module.
The modules $\UGHL{\Lambda}$ and $\FGHL{\Lambda}$ are universal and final, respectively, in the following sense.
\begin{propn}\label{UniFinalPropPropn}
Let $N$ denote a weight $\g$-module with generative highest layer $\Gamma$, and let $\Lambda$ denote a weight $\Heis$-module.  Then:
\begin{enumerate}
\item Any epimorphism of $\Heis$-modules $\Lambda \EpiArrow \Gamma$ extends uniquely to an epimorphism of $\g$-modules
$$ \UGHL{\Lambda} \EpiArrow N;$$
\item Any epimorphism of $\Heis$-modules $\Gamma \EpiArrow \Lambda$ extends uniquely to an epimorphism of $\g$-modules
$$ N \EpiArrow \FGHL{\Lambda}.$$
\end{enumerate}
\end{propn}
For any function $\varphi: \Z \to \K$, write $\FGHL{\varphi}$ for the $\g$-module $\FGHL{\PolyMod{\varphi}}$.
\begin{propn}\cite{Chari}
For any two functions $\varphi, \psi: \Z \to \K$, the $\g$-modules $\FGHL{\varphi}$ and $\FGHL{\psi}$ are isomorphic if and only if there exists a non-zero $\lambda \in \K$ such that
$$ \varphi (k) = \lambda^k \hspace{0.1em} \psi(k), \quad k \in \Z.$$
\end{propn}
\begin{propn}\label{ISQhaveGHLProposition}
Let $N$ denote an irreducible subquotient of the $\g$-module $\VermaQuotient{0}$.  Then:
\begin{enumerate}
\item \label{ISQGHL1} There exists a non-negative integer $n$ such that $N_n$ is the generative highest layer of $N$;
\item \label{ISQGHL2} The layer $N_n$ is an irreducible subquotient of the $\Heis$-module $\Level{n}$;
\item \label{ISQGHL3} The $\g$-modules $N$ and $\FGHL{N_n}$ are canonically isomorphic. 
\end{enumerate}
\end{propn}
\begin{proof}
As $N$ is a subquotient of $\VermaQuotient{0}$, it has layer decomposition
$$ N = \oplus_{n \geqslant 0}{N_n}.$$
Let $n \geqslant 0$ be minimal such that $N_n \ne 0$.
Then $\E \cdot N_n = 0$ since $N_{n-1} = 0$, and $N = \UEA{\g}N_n$ since $N$ is irreducible.
Hence $N_n$ is the generative highest layer of $N$, proving part \ref{ISQGHL1}.

Now let $P' \subset P$ be $\g$-submodules of $\VermaQuotient{0}$ such that $N = P / {P'}$.
Then $N_l = P_l / {P'_l}$, for all $l \geqslant 0$, and in particular $N_n$ is a subquotient of the $\Heis$-module $\Level{n}$.
Suppose that $N_n$ is a reducible $\Heis$-module, and let $\Lambda \subset N_n$ denote a proper submodule.
By the Poincar\'e-Birkhoff-Witt Theorem,
\begin{eqnarray*}
\UEA{\g} \Lambda &=& \UEA{\F} \otimes \UEA{\Heis} \otimes \UEA{\E} \Lambda \\
	&=& \UEA{\F} \Lambda,
\end{eqnarray*}
and so $N_n \not \subset \UEA{\g} \Lambda$.
Hence $\UEA{\g}{\Lambda}$ is a proper $\g$-submodule of $N$, contrary to hypothesis.
Therefore $N_n$ is an irreducible $\Heis$-module, and so part \ref{ISQGHL2} is proven.

The $\g$-module $N$ has generative highest layer $N_n$, and so by Proposition \ref{UniFinalPropPropn} there is an epimorphism of $\g$-modules
$$ N \EpiArrow \FGHL{N_n}.$$
As $N$ is irreducible, this map is an isomorphism, and so part \ref{ISQGHL3} is proven.
\end{proof}
\begin{theorem}\label{VQISQClassnTheorem}
For any $\varphi \in \ExpFunc$, the $\g$-module $\FGHL{\varphi}$ is an irreducible subquotient of $\VermaQuotient{0}$.
Moreover, if $\K$ is algebraically closed, then any irreducible subquotient of $\VermaQuotient{0}$ is of the form $\FGHL{\varphi}$ for some $\varphi \in \ExpFunc$.
\end{theorem}
\begin{proof}
Let $n$ be a positive integer, and let $\varphi \in \ExpFunc_n$.
Theorem \ref{ExistenceSV} guarantees the existence of a non-zero singular vector $v \in \Level{n}$.
The $\Heis$-module $\UEA{\Heis}v$ contains only singular vectors, and is isomorphic to $\Level{n}$, by Corollary \ref{UHvIsoMr}.
Let $P = \UEA{\g}v$.  By the Poincar\'e-Birkhoff-Witt Theorem, 
\begin{eqnarray*}
\UEA{\g} \cdot v &=& \UEA{\F} \otimes \UEA{\Heis} \otimes \UEA{\E} \cdot v \\
	&=& \UEA{\F} \otimes \UEA{\Heis} \cdot v
\end{eqnarray*}
as vector spaces.
Therefore $P = \oplus_{l \geqslant n}{P_l}$ and $P$ is a $\g$-submodule of $\VermaQuotient{0}$ with generative highest layer $P_n = \UEA{\Heis}v \cong \Level{n}$.
Hence, by Theorem \ref{LevelImageClassification}, there is an epimorphism of $\Heis$-modules
\begin{equation*}
	P_n \EpiArrow \PolyMod{\varphi},
\end{equation*}
which by Proposition \ref{UniFinalPropPropn} extends to an epimorphism of $\g$-modules
\begin{equation*}
P \EpiArrow \FGHL{\varphi}.
\end{equation*}
Thus $\FGHL{\varphi}$ is an irreducible subquotient of $\VermaQuotient{0}$.

Now suppose that $\K$ is algebraically closed, and that $N$ is an irreducible subquotient of $\VermaQuotient{0}$.
Then by Proposition \ref{ISQhaveGHLProposition}, there exists $n \ge 0$ such that 
\begin{enumerate} \item $N_n$ is the generative highest layer of $N$;
\item $N_n$ is an irreducible subquotient of the $\Heis$-module $\Level{n}$;
\item $\FGHL{N_n} \cong N$. 
\end{enumerate}
By Theorem \ref{LevelImageClassification}, there exists a $\varphi \in \ExpFunc_n$ such that $N_n \cong \PolyMod{\varphi}$ as $\Heis$-modules, and so
$$ N \cong \FGHL{N_n} \cong \FGHL{\varphi},$$
which completes the proof of the Theorem.
\end{proof}

The following Corollary is immediate from \cite{BilligZhao} and the observation that $\ExpFunc$ is a family of exp-polynomial functions.
\begin{corollary}
Suppose that $\K$ is algebraically closed.  
Then the components of any irreducible subquotient of $\VermaQuotient{0}$ have finite dimension.
\end{corollary}
\end{section}

\begin{section}{Acknowledgements}
The author is indebted to Vyacheslav Futorny and Alexander Molev for their careful guidance, and would like to thank Yuly Billig, Gus Lehrer and Paulo Agonozzi Martin for their assistance.

This work was completed as part of a cotutelle Ph.~D.~ programme at the Universidade de S\~ao Paulo and the University of Sydney.
\end{section}


\end{document}